\documentclass[reqno]{amsgok}
\usepackage{epsfig}

\newtheorem{thm}{Theorem}[section]
\newtheorem{cor}[thm]{Corollary}

\newtheorem{prop}[thm]{Proposition}

\theoremstyle{definition}

\theoremstyle{remark}
\newtheorem{rem}{Remark}[section]

\newcommand{\fb}{\fbox} 
\let\fb=\relax

\begin{document}

\title{Exotic structures and adjunction inequality}
\author[Akbulut \& Matveyev]{Selman Akbulut 
\footnote {The first author was partially supported by an NSF grant DMS-9626204} and Rostislav Matveyev}
\address{Michigan State University, Dept of Mathematics, E. Lansing MI 48824, USA}
\email{akbulut@math.msu.edu \\ matveyev@math.msu.edu}

\volume{5}
\fpage{}
\lpage{}



\maketitle


\section {Introduction}

Here we want to reprove and strengthen some old difficult 
theorems of the theory of $4$-manifolds by the aid of recently proven modern tools. 
One of the important recent results of smooth 
$4$-manifolds is Eliashberg's topological characterisation of 
compact Stein manifolds, that is complex manifolds which admit strictly
plurisubharmonic Morse function:

\begin{thm}(\cite{E})
Let $X=B^4\cup(\text{1-handles})\cup(\text{2-handles})$ 
be four-dimensional handlebody with one 
0-handle and no $3$- or $4$-handles. 
Then 
\begin{itemize}
\item The standard symplectic structure on $B^4$ can 
      be extended over $1$-handles so that 
      manifold $X_1=B^4\cup(\text{1-handles})$ 
      is a compact Stein domain.
\item If each $2$-handle is attached to 
      $\partial X_1$ along a Legendrian 
      knot  with framing one 
      less then Thurston-Bennequin framing 
      of this knot, then the symplectic form
      and complex structure 
      on $X_1$ can be extended over 
      $2$-handles to a symplectic form on $X$,
      so that $X$  become compact Stein domain.
\end{itemize}
\end{thm}

P. Lisca and G. Mati\'c showed that compact Stein manifolds imbed naturally  into
K\"ahler surfaces:

\begin{thm}(\cite{LM})
Interior of every compact Stein manifold $X$ can be holomorphically embedded as a domain 
into a minimal K\"ahler surface $S$ with ample canonical bundle and 
$\;b_{2}^{+}(S) >1$, such that the induced 
K\"ahler form on $X$ agrees with the original symplectic form on $X$.
\end{thm}

\begin{thm}
A minimal K\"ahler surface $X$, with $ \; b_{2}^{+}(X) >1 $ and an ample 
canonical bundle,
can not contain a smoothly embedded $2$-sphere $\Sigma \subset X$ 
with $\Sigma . \Sigma \geq -1$ .
\end{thm}

This theorem roughly follows from the fact that K\"ahler surfaces have nonzero 
Seiberg-Witten invariants (see \cite {B} and \cite {MF} in 
case $ \Sigma . \Sigma =-1$, and 
 \cite{FS} in case $\Sigma . \Sigma \geq 0$). Also embedded surfaces
in K\"ahler 
manifolds satisfy the  
adjunction inequality:

\begin{thm} (\cite{KM1}, \cite{MST})
Let $X$ be a closed smooth $4$-manifold with  $ b_{2}^{+}(X) >1$, with
a nonzero Seiberg Witten invariant (e.g. $X$ K\"ahler surface) corresponding 
to the line bundle $L\to X$. Let 
$\Sigma \subset X$ be a compact oriented embedded surface with 
$\Sigma.\Sigma \geq 0$ and $\Sigma $ is not diffeomorphic to a sphere, then
$$ 2g(\Sigma )-2 \geq \Sigma. \Sigma + |c_{1}(L).\Sigma | $$  
\end{thm}

After discussing the basics, we will show some $4$-manifold theorems 
can be obtained as easy corollaries of these three basic theorems.

\section[]{Definitions}

Let us recall some basic facts (see, for example \cite{G}, \cite{E}). 
Any compact Stein domain $X$ 
induces a contact structure $\xi$ on the boundary
$3$-manifold $Y=\partial X$. It is isomorphic to the restriction 
of the dual $ K^{*}$ of the canonical 
line bundle $ K\to X$. Furthermore if  $\alpha $ is an oriented
Legendrian knot in $ Y $ (a knot whose tangents lie in the contact planes)
 bounding an oriented surface $ F $ in $ X $, then the
``rotational number'' is defined to be the relative Chern
class  $ rot(\alpha ,F)= c_{1}({K}^{*},v)$ of the induced  $2$-plane 
bundle  ${K}^{*} \to F$  with respect to the tangent vector field 
$v$ of $\alpha $ (i.e. the obstruction to extending $v$ to a section of 
$K^{*}$ over $F$). Contact structure $\xi$ also induces the so called
Thurston-Bennequin framing $tb(\alpha)$ by taking normal vectorfield $u$
along $\alpha$, such that $u\in \xi$.

The simplest example of a compact Stein domain is $ B^{4}\subset {\bf C}^{2} $ with
the induced symplectic structure. One can choose coordinates in
$ {\bf R}^{3}\subset  S^{3} = \partial B^{4} $, so that the induced 
contact structure $\xi_{0}$ on ${\bf R}^{3}$ is the
kernel of the form $ \lambda_{0}=dz+xdy $. By Theorem 1.1 this
symplectic structure 
extends across $1$-handles attached to $B^{4}$, such that resulting
manifold is comact Stein domain;
we draw the attaching balls of each $1$-handle on the line 
$\{x=0,z=\mbox{constant}\}$. Any link $L$ in $S^{3}$ 
can be isotoped to a Legendrian link, we can achieve this by first isotoping $L$
so that all crossings in the front projection of $L$ 
are left handed as in the first picture of Figure 1.
For example, we can turn a right handed crossing to a left handed 
crossing by the local isotopy as in the second picture of Figure 1.

\begin{figure}[htb]
\fb{\epsfig{figure=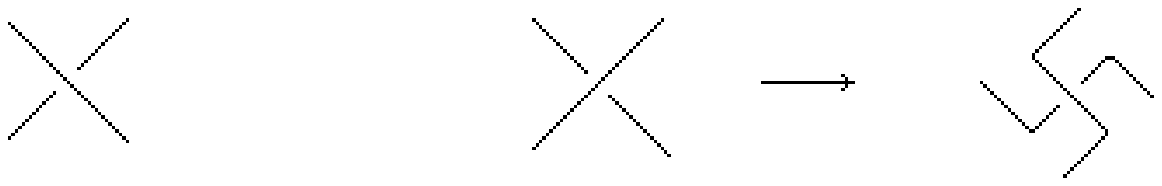, width=100mm}}
\caption{}
\end{figure}

Projection of a Legendrian link to the $yz$-plane (front projection) has
intersections only coming from left handed crossing, no vertical tangencies
and all minima and maxima in $y$-direction are cusps (as the projection
of the knot in Figure 3 to the $yz$-plane). Moreover, every projection with 
these properties is a projection of some Legendrian link.
Rotational number of an oriented Legendrian knot $\alpha $ in $S^{3}$ does not
depend on the surface $F\subset S^{3}$ it bounds. Invariants $rot(K)$ and $tb(K)$ 
can be calculated by:
$$ rot(\alpha) = 1/2 (\mbox{ Number of ``downward"  cusps } - \mbox
{ Number of ``upward" cusps })$$
 $$ tb(\alpha)=bb(\alpha)- c(\alpha )$$ 
where $bb(\alpha)$ is the blackboard ($yz$-plane) framing of the projection 
of $\alpha $, $c(\alpha )$ is the number of right cusps, and 
``downward" and ``upward" cusps are calculated in the obvious way.  
 Because the local 
isotopy  in Figure 1 introduces one right cusp (hence a $-1$ contribution 
to calculation of tb$(K)$) we can incorporate this to ``self crossing number" 
calculation by reading the crossings numbers in a modified way  as in in Figure 2.  

\begin{figure}[htb]
\fb{\epsfig{figure=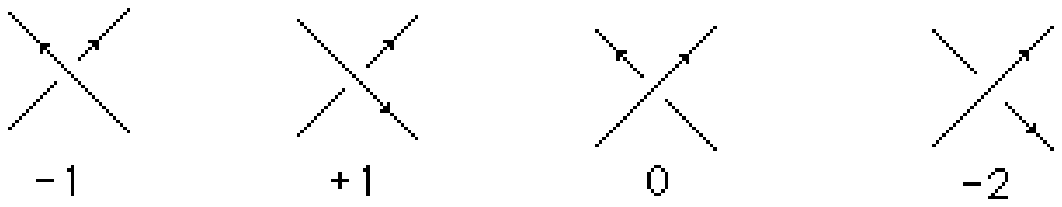, width=100mm}}
\caption{}
\end{figure} 

For example in the knot of Figure 3 we can calculate tb$(K)=5-4=1 $ and 
rot$(K)=0 $. A useful corollary to these theorems is the following generalization of 
the Bennequin inequality: 
 
\begin{figure}[t]
\fb{\epsfig{figure=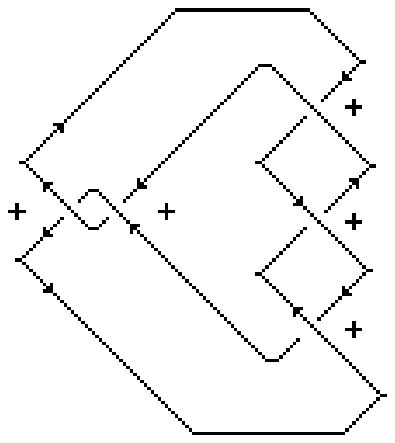, width=50mm}}
\caption{}
\end{figure}

\begin{cor}
Let $X$ be a compact Stein domain, $F\subset X$ be 
a two-dimensional submanifold of $X$, such that 
$\alpha=\partial F\subset\partial X$ is  Legendrian 
with respect to induced contact structure and $f$ is
framing on $\alpha$ induced by a trivialization of
the normal bundle of $F$ in $X$, then 
$$-\chi(F) \geq [\mbox{tb}(\alpha)-f]+|\mbox{rot}(\alpha)|$$
where  tb$(\alpha)$ and rot$(\alpha) $ are the 
Thurston-Bennequin framing and 
rotational number of $\alpha $  
\end{cor}

To prove this we  attach a $2$- handle to $X$ along 
$\alpha $ with the framing $tb(\alpha)-1 $, and apply  Theorems 1.1, 1.2, and 1.3
to the resulting manifold 
and closed surface $F'=F\cup_{\partial} D$, where $D$ is the core $2$-disc 
of the $2$-handle.

\section[]{Applications}
\subsection{}
\begin{thm}(\cite{A1})
Let $W$ be the contractible manifold of Figure 4. Let
$f:\partial W \rightarrow \partial W $, be the involution induced by an 
involution of $S^{3}$ (as described in \cite{A1}) with 
$f(\gamma)=\gamma'$, where 
$\gamma$ and $\gamma'$ are circles in $\partial W$ 
shown on Figure 4.
Then $f:\partial W \rightarrow \partial W $ does not extend to a diffeomorphism
$F: W\rightarrow W$ (but it extends to a homeomorphism) 
\end{thm}

\begin{figure}[b]
\fb{\epsfig{figure=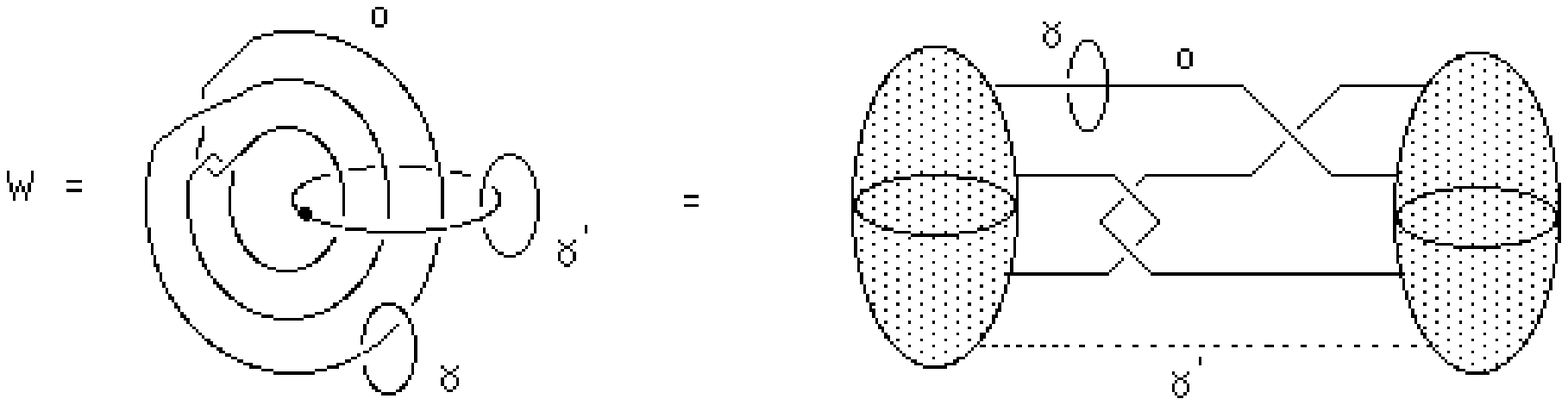, width=110mm}}
\caption{}
\end{figure} 

\begin{proof}
 By applying Theorem 1.1 to manifold $W$ 
(the second picture of $W$ in Figure 4) we see that $W$ 
is compact Stein. Also since $\gamma $ is slice in $W$ if $f$ extended to 
a diffeomorphism $F: W\to W$, then
$\gamma '$ would be slice also. But this contradicts the inequality of Theorem 2.1 
(here $F=D^{2}$, $f=0$, and $tb(\gamma ')=0$) 
\end{proof}

\subsection{}
\begin{thm}(\cite{A2})
Let $Q_{1}$ and $Q_{2}$ be the manifolds obtained by attaching  $2$-handles to 
$B^{4}$ along the knots $K_{1}$ and $K_{2}$ with $-1$ framings as in Figure 5, then
$Q_{1}$ and $Q_{2}$ are homeomorphic but not diffeomorphic to each other, 
even interiors are not diffeomorphic to each other.
\end{thm}

\begin{figure}[htb]
\fb{\epsfig{figure=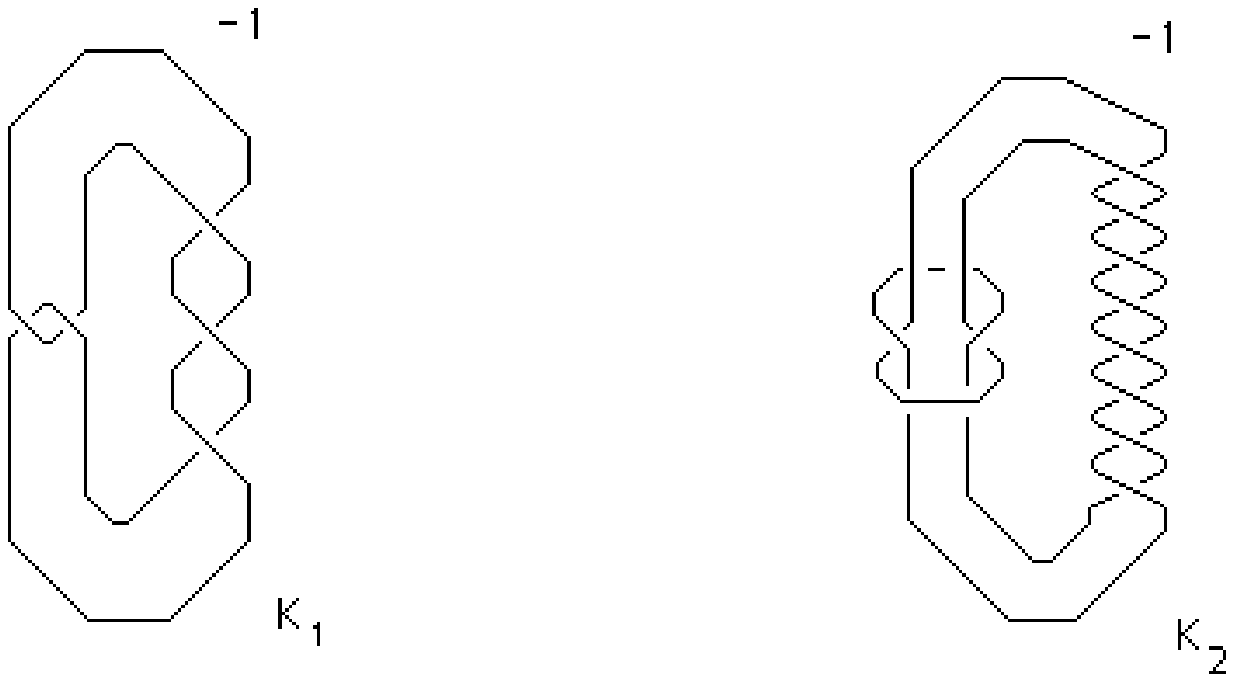, width=100mm}}
\caption{}
\end{figure}

\begin{proof}

By Theorem 1.1 manifold $Q_{1}$ is compact Stein domain, by Theorem
1.2 it embeds as a domain
in a minimal K\"ahler surface $S$ with $b_{2}^{+}(S)>1$. If $Q_{1}$ were diffeomorphic 
to $Q_{2}$, the generator of $H_{2}(Q_{1};{\bf Z})$ would be represented by a 
smooth embedded sphere (since $K_{2}$ is a slice knot the generator 
of $H_{2}(Q_{1};{\bf Z})$ represented by a smooth sphere) with $-1$ self intersection,
violating Theorem 1.3. 
\end{proof}

\subsection{}
Let $K\subset S^{3}$ be a Legendrian knot, and $K_{0}'$ be a $0$-push off of $K$
(the zero framing is the framing induced from the normal vector field of $K$ in
the oriented surface in $S^{3}$ bounding $K$ ).

\begin{prop} 
We can move $K_{0}'$ to a Legendrian knot 
$K_{0}$  by an isotopy which fixes $K$, such that
$tb(K_{0})=-|tb(K)|$
\end{prop}

\begin{proof}

$tb(K)=bb(K)-c(K)$. If $ bb(K)\leq 0$ then $K_{0}'$ is just the blackboard push-off of 
$K$ with $-2bb(K)$ right half twist, but this is a projection of a Legendrian link
and $bb(K_{0})=bb(K)$ and $c(K_{0})=c(K)$ hence $tb(K_{0})=tb(K)$

\begin{figure}[htb]
\fb{\epsfig{figure=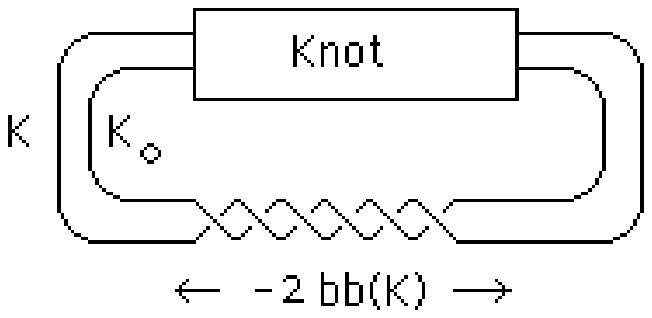, width=50mm}}
\caption{}
\end{figure} 

 If $tb(K)\leq 0$ and $ bb(K) > 0$, then $K_{0}'$ is the blackboard push-off  with 
$2bb(K)$ left handed half twist. Knot $K$ has $2c(K)$ cusps, and $2c(K)\geq 2bb(K)$. 
We can produce 
a Legendrian picture by placing all half twists near cusps as in Figure 7.
Thus $tb(K_{0})=tb(K)$.

\begin{figure}[htb]
\fb{\epsfig{figure=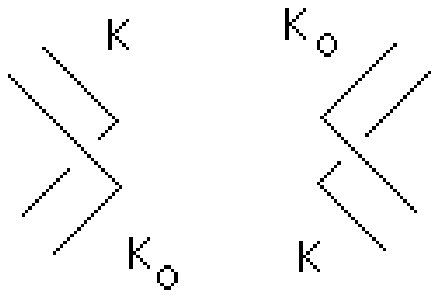, width=50mm}}
\caption{}
\end{figure}

 If $tb(K)>0$, then $bb(K)>0$ and again $K_{0}'$ is the blackboard push-off 
with $2bb(K)$ left half twist. By above trick of placing
half left twist near cusps, we can get rid of $2c(K)$ left twist, and we are left with
$ 2bb(K)-2c(K)$ left twist each of which contributes  $-1$  to $tb(K_{0})$
 as in Figure 8.
Hence $tb(K_{0})=tb(K)-[\;2bb(K)-2c(K)\;]=tb(K)-2tb(K)=-tb(K)$ 

\begin{figure}[b]
\fb{\epsfig{figure=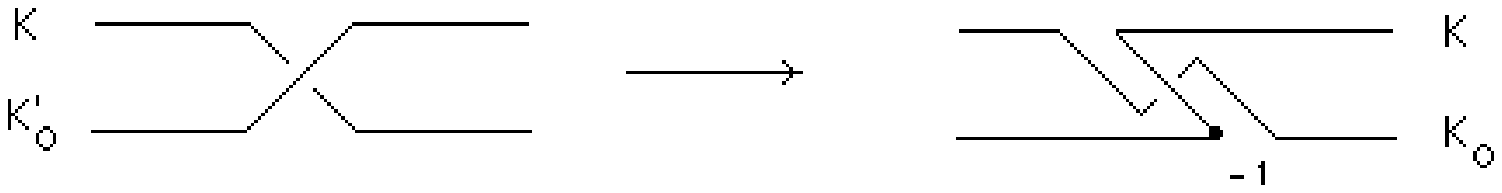, width=120mm}}
\caption{}
\end{figure}

\end{proof}

By Proposition 3.3 we can produce many $0$-push offs $K_{0}^{i}$, $i=1,2,..k$
such that  $tb(K_{0}^{i})=-|tb(K)|$. So if $tb(K)\leq 0$ then all $K_{0}^{i}$ 
have same $tb$, if $tb(K)>0$ 
then all $K_{0}^{i}$ but the original knot $K$ has the same $tb$. 

\begin{thm} 
If $K\subset S^{3}$ a Legendrian knot with $tb(K)\geq 0$ , then 
all iterated positive Whitehead doubles of $Wh_{n}(K)$ are not slice. In fact if 
$Q_{n}^{r}(K)$ is the
manifold obtained by attaching a $2$-handle to $B^{4}$ along $Wh_{n}(K)$ 
 with a framing $-1\leq r\leq 0$, then
there is no smoothly embedded $2$-sphere
in $Q_{n}^{r}(K)$ representing the generator of 
$H_{2}(Q_{n}^{r}(K);{\bf Z})$. 
\end{thm}

\begin{proof}
The first positive Whitehead double is obtained by connecting $K$ and $K_{0}$ by a 
left handed cusp as in Figure 9 which contributes $+1$ to $tb$, hence
$$tb(Wh(K))=tb(K)+tb(K_{0})+1=tb(K)-tb(K) +1=1$$
So by iteration we get $tb(Wh_{n}(K))=1$. But Corollary 2.1 says that any slice knot
$L$ must have $tb(L)\leq -1$.

\begin{figure}[htb]
\fb{\epsfig{figure=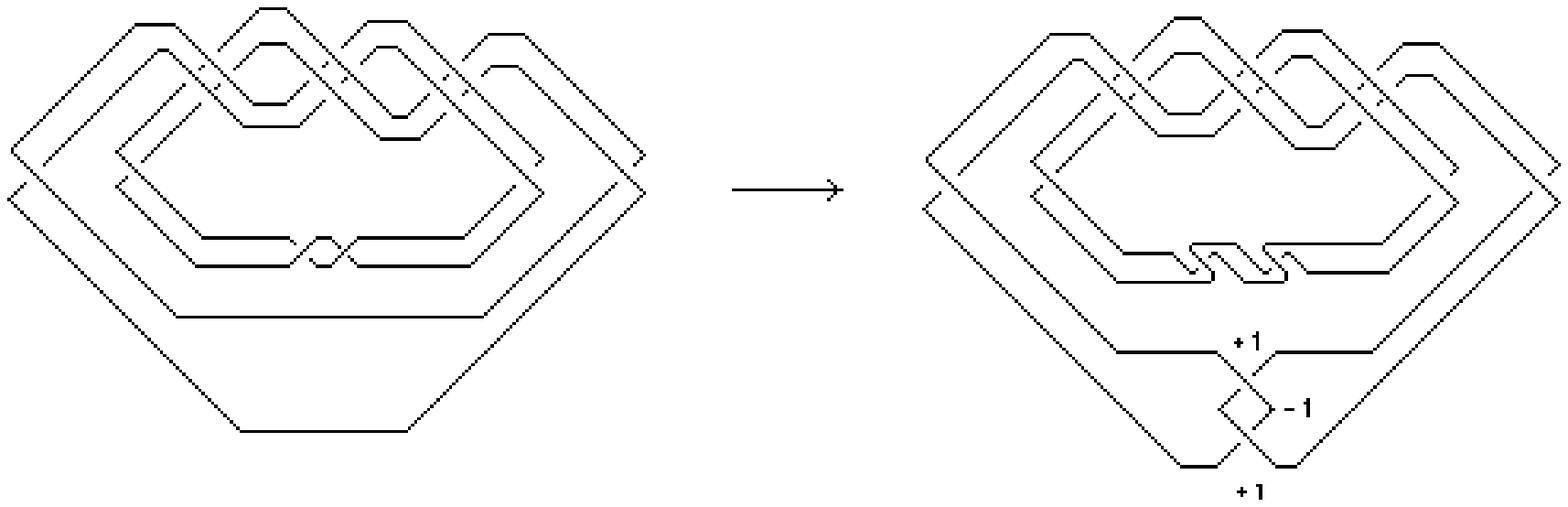, width = 133mm}}
\caption{}
\end{figure}

For the second part, observe that since $ r\leq 0 = tb(Wh_{n}(K))-1 $, by Theorem 1.1   
$Q_{n}^{r}(K)$ is compact Stein. By Theorem 1.2 we can imbed $Q_{n}^{r}(K)$ into 
K\"ahler surface $S$. Then by Theorem 1.3  we can not have a smoothly
embedded $2$-sphere $\Sigma\subset Q_{n}^{r}(K)\subset S$
 representing $H_{2}(Q_{n}^{r}(K);{\bf Z})$, since $\Sigma . \Sigma =r \geq -1$.
\end{proof}

\begin{rem}
 L.Rudolf has previously shown that $Wh_{n}(K)$ are not slice 
if $tb(K)\geq 0$ (\cite {R}).
\end{rem}

\begin{rem}
Proof of Proposition 3.2 gives: If $K_{r}'$ is the $r$-framing push off 
of a knot $K$, then $K_{r}'$ can be isotoped to a Legendrian knot $K_{r}$ fixing
$K$ with $(tb(K_{r})-r)=-|tb(K)-r|$
\end{rem}

\begin{rem}
If $\xi$ is the contact structure 
on  $\Sigma =\partial W$  induced by the symplectic structure on $W$ of Theorem 3.1, 
and $f^{*}(\xi)$  is the ``pull-back" contact structure on $\Sigma $. 
Then it follows that the contact structures $\xi$ and $f^{*}(\xi)$ are homotopic 
through $2$-plane fields but not isotopic through contact structures (\cite {AM}). 
\end{rem}

\end{document}